\documentclass[12pt]{article}
\usepackage{amssymb,amsmath,graphicx,amsthm}
\usepackage[T1]{fontenc}
\usepackage{ae,aecompl}
\newtheorem{theorem}{Theorem}
\newtheorem{lemma}[theorem]{Lemma}

\makeatletter
\renewcommand{\@makefnmark}{\hbox{\mathsurround=0pt}$^{\mbox{\dag}}$}
\makeatother

\begin{document}

{\Large\noindent Through the Looking-Glass of the Grazing Bifurcation:\\[3pt] Part I - Theoretical Framework}

\bigskip

\setcounter{footnote}{1}
{\noindent\large James Ing${}^1$, Sergey Kryzhevich${}^2$\footnote{Corresponding author.\\\hspace*{5.5mm}Email address:
kryzhevitz@rambler.ru }~, Marian Wiercigroch${}^1$}

\bigskip


{\noindent\normalsize ${}^1$ Centre for Applied Dynamics Research, School of Engineering, University of Aberdeen, \\ Kings College Aberdeen AB24 3UE, Scotland, UK;\\
${}^2$ Faculty of Mathematics and Mechanics, Saint-Petersburg State University,\\
28, Universitetskiy pr., Peterhof, Saint-Petersburg, 198503, Russia,\\
University of Aveiro, Department of Mathematics, 3810$-$193, Aveiro, Portugal}




\bigskip

\begin{table}[h!]
\vspace*{-5mm}
\doublerulesep 0.05pt
\tabcolsep 7.8mm
\vspace*{2mm}
\setlength{\tabcolsep}{7.5pt}
\hspace*{-2.5mm}\begin{tabular*}{16.5cm}{r|||||l}
\multicolumn{2}{l}{\rule[-6pt]{16.5cm}{.01pt}}\\
\parbox[t]{6cm}{\small
\vspace*{.5mm}
\hfill {\bf Submission Info}\par
\vspace*{2mm}
\hfill Communicated by Referees\par
\hfill Received DAY MON YEAR \par
\hfill Accepted DAY MON YEAR\par
\hfill Available online DAY MON YEAR\par
\noindent\rule[-2pt]{6.3cm}{.1pt}\par
\vspace*{2mm}
\hfill {\bf Keywords}\par
\vspace*{2mm}
\hfill Grazing \par
\hfill Homoclinic point \par
\hfill Structural stability \par
\hfill Models of impact}
&
\parbox[t]{9.85cm}{
\vspace*{.5mm}
{\normalsize\bf Abstract}\par
\renewcommand{\baselinestretch}{.8}
\normalsize \vspace*{2mm} {\small It is well-known for vibro-impact systems that the existence of a periodic solution with a low-velocity impact (so-called grazing) may yield complex behavior of the solutions. In this paper we show that unstable periodic motions which pass near the delimiter without touching it may give birth to chaotic behavior of nearby solutions. We demonstrate that the number of impacts over a period of forcing varies in a small neighborhood of such periodic motions. This allows us to use the technique of symbolic dynamics. It is shown that chaos may be observed in a two-sided neighborhood of grazing and this bifurcation manifests at least two distinct ways to a complex behavior. In the second part of the paper we study the robustness of this phenomenon. Particularly, we show that the same effect can be observed in "soft"\ models of impacts.}
\par
\par
\par
\hfill{\scriptsize\copyright 2012 L\&H Scientific Publishing, LLC. All rights reserved.}}\\[-4mm]
&\\
\multicolumn{2}{l}{\rule[15pt]{16.5cm}{.01pt}}\\\end{tabular*}
\vspace*{-7mm}
\end{table}


\renewcommand{\baselinestretch}{1}
\normalsize

\bigskip

\noindent\textbf{1. Introduction}

\bigskip

\noindent Vibro-impact systems appear in different mechanical problems (modeling of impact dampers, clock mechanisms, immersion of constructions, etc.). Their properties show a lot of resemblence to classical nonlinear systems. Particularly, chaotic dynamics is possible [1--21].

The so-called grazing bifurcation was first described by A.\, Nordmark \cite{iwk6,iwk15}. The critical value of the parameter corresponds to a zero velocity impact of a periodic solution. It was demonstrated that this bifurcation implies non-smooth behavior of solutions, instability of the periodic solution in the parametric neighborhood of grazing and, in additional assumptions, chaotic dynamics. It was shown theoretically, how periodic solutions corresponding to different values of impacts $(n_1,n_2,\ldots, n_m)$ over successive periods of forcing may appear \cite{iwk16}.

However, parametric neighborhoods of the so-called continuous grazing \cite{iwk15,iwk16} consist of two parts. One of them corresponds to periodic solutions with a low-velocity impact, another one corresponds to periodic motions that pass close to the delimiter but do not touch it. In our paper we consider the second case assuming that the mentioned periodic motions are unstable.

We consider the system studied in \cite{iwk13}: the variable $x\in [-1,1]$ satisfies the equation
\begin{equation}\label{eiwk1.1}
\ddot x + 2\beta \dot x-x=f(t)
\end{equation}
over free flight intervals. Impacts, given by the equalities
$$v_+=-rv_-\quad \mbox{where} \quad \dot v_\pm=\dot x (t\pm 0)$$
take place at $|x|=1$. It was assumed that the right hand side of the equation \eqref{eiwk1.1} can be represented as a combination of Dirac functions,
$$f(t)=\gamma\sum_{k=-\infty}^{\infty} (\delta(t-kT-a)-\delta(t-kT)),$$
where $a\in (0,T)$. For the considered system the authors proved the existence of so-called non-classical bifurcations, corresponding to coincidence of the impulse action which takes place for $t=a+kT$ and $t=kT$ and the impact. These bifurcations are accompanied by the disappearance of the periodic solution and, sometimes, by appearance of a "strange attractor".

The main aim of our article is to show that chaotic dynamics is possible in a two-sided neighborhood of grazing. Namely, it can be caused by the presence of a unstable periodic motion passing near the delimiter. An experimental illustration of this is given in \cite{iwk14}.

We start by considering an impulse model of impacting system where interaction between the moving particle and the delimiter is instantaneous. However, from a physical point of view the so-called "soft"\ model \cite{iwk2,iwk9,iwk10,iwk14} where the delimiter is considered as a very stiff spring, is more relevant. We prove that any hyperbolic invariant set that appears in the impulse model persists if we pass to a "soft model" with a sufficiently stiff spring.

The rest of the paper is organized as follows. In Section 2 we give a description of the impulse model of impacts. In Section 3 the grazing family of periodic solutions is considered. Section 4 is devoted to description of a set of initial conditions corresponding to solutions with a zero velocity impact. In Section 5 the main result on the existence of chaotic invariant set is formulated. In Section 6 asymptotic estimates for Jacobian matrices corresponding to near-grazing solutions are given. The existence of a transverse homoclinic point is proved at Secition 7, and in the next section the corresponding symbolic dynamics is described. A simple example of a single degree-of-freedom system is considered in Section 9. Also we prove in Section 10 that an infinitely stiff delimiter may be replaced with a sufficiently stiff spring. The conclusion is given in Section 11.

Later on we use the following formalism: indices of parameters of successive impacts of a fixed solution (phases, velocities, etc.) are denoted by superscripts in order to distinguish them from ones of coordinates of vectors always denoted by subscripts. Denote by ${\mathop {\rm col}\nolimits} (a_1,\ldots,a_m)$ the column vector, consisting of elements $a_1,\ldots, a_m$. Here $a_k$ themselves may be vectors. Also, we use the notation $\bar a={\mathop {\rm col}} (a_2,\ldots a_n)$.

\bigskip

\noindent\textbf{2. Impulse model}

\bigskip

\noindent We study the motion of a point mass, described by system of second order differential equations of the general form and impact conditions of impulse type.

Consider a segment $J=[\mu_-,\mu_+]$ and a $C^1$ smooth function
$${ f}(t,{ z},\mu):{\mathbb R}^{2n+1}\times J\to {\mathbb R}^n.$$
Here $\mu$ is a scalar parameter. For example, this can be the amplitude of free periodic oscillations of the considered system or the coefficient of restitution. Suppose that
$${ f}(t,{ z},\mu)\equiv { f}(t+T,{ z},\mu).$$
Let
$$\begin{array}{c}
{ z}={\mathop {\rm col}\nolimits} ({ z}_{1}, \ldots, { z}_{n}); \quad
{ z}_k={\mathop {\rm col}\nolimits} (x_k, y_k),\quad k=1,\ldots, n; \quad {\bar { z}}={\mathop {\rm col}\nolimits}({ z}_{2},\ldots, { z}_{n});\\[5pt]
\quad { f}={\mathop {\rm col}\nolimits}(f_{1},\ldots, f_{n});\quad
{ x}={\mathop {\rm col}\nolimits}(x_{1}, \ldots, x_{n}); \quad { y}={\mathop {\rm col}\nolimits}(y_{1}, \ldots, y_{n});\\[5pt] \quad {\bar{ x}}={\mathop {\rm col}\nolimits}(x_{2}, \ldots, x_{n}),\quad
{\bar { y}}={\mathop {\rm col}\nolimits}(y_{2}, \ldots, y_{n}).\end{array}$$
Consider the system
\begin{equation}\label{eiwk2.1}
\dot x_k=y_k;\qquad \dot y_k=f_k(t,{ z},\mu), \quad k=1,\ldots, n.
\end{equation}

Let the coefficient of restitution $r=r(\mu)\in (0,1]$ be a $C^1$~-- smooth function. Suppose that Eq. \eqref{eiwk2.1} is defined for ${ z}\in\Lambda=[0,+\infty)\times {\mathbb R}^{2n-1}$ and the following impulse type impact conditions take place as the component $x_1$ vanishes.

\textbf{Condition 1.} \
\begin{enumerate}
\item If $x_1(t_0)=0$ then ${ x}(t_0+0)={ x}(t_0-0)$,
\begin{equation}\label{eiwk2.2}
y_1(t_0+0)=-r(\mu)y_1(t_0-0),\qquad {\bar y}(t_0+0)=\bar { y}(t_0-0).
\end{equation}
\item $x_1(t)\geqslant 0$ for all $t$ where ${ z}(t)$ is well-defined.
\end{enumerate}

Consider the vibro-impact system

\begin{equation}\label{eiwk2.3}\begin{array}{c}
\dot x=y;\qquad \dot y=f(t,{ z},\mu);\\
\mbox{Condition 1 is applied if } x_1=0.
\end{array}\end{equation}

We give two definitions of solutions of vibro-impact systems.

\textbf{Definition 1.} The function ${ z}(t)=\mbox{\rm col}(x_1(t),y_1(t),\ldots,x_n(t),y_n(t))$ is a solution of Eq. \eqref{eiwk2.3} \emph{with a finite number of impacts over the interval $(a,b)$}, if there exists a finite number of instants $a=t_0<t_1<\ldots<t_N<t_{N+1}=b$ such that the following conditions are satisfied.
\begin{enumerate}
\item All components of the solution ${ z}(t)$, except $y_1(t)$, are continuous while $t\in (a,b)$. The discontinuity set of the function $y_1(t)$ is a subset of $\{t_k\}$.
\item The function $x_1(t)$ is non-negative on $(a,b)$. The set $\{t_1,\ldots,t_N\}$ is the set of all zeros for this function.
\item For any $k=1,\ldots,N$
$$y_1(t_k+0)=-r(\mu)\,y_1(t_k-0).$$
\item The function ${ z}(t)$ is a solution of \eqref{eiwk2.1} on every interval $(t_k,t_{k+1})$.
\end{enumerate}

For completeness of the mathematical model, we define a solution with an infinite number of impacts.

\textbf{Definition 2.} The function ${ z}(t)=\mbox{\rm col}(x_1(t),y_1(t),\ldots,x_n(t),y_n(t))$ is a solution of the vibro-impact system \eqref{eiwk2.3} over the interval $(a,b)$ if there exists a disjoint splitting $(a,b)={\cal I}_+\bigcup {\cal I}_0 \bigcup {\cal I}_-$ with the following properties.
\begin{enumerate}
\item The set ${\cal I}_+=\{t\in (a,b):x_1(t)>0\}$
corresponding to free flight motions is open. The set ${\cal I}_-=$
$$ \{t\in (a,b): x_1(t)=y_1(t)=0, f_1(t,0,0,x_2(t),y_2(t),\ldots,x_n(t),y_n(t),\mu)\leqslant 0\}$$
is closed.
\item All the components of the solution ${ z}(t)$, except $y_1(t)$, are continuous over $(a,b)$, the discontinuity set of the function $y_1(t)$ is a subset of $I_0$.
\item The function $x_1(t)$ is non-negative on $(a,b)$. The set ${\cal I}_0\bigcup {\cal I}_-$ is the set of zeros of this function.
\item The function ${ z}(t)$ is a solution of \eqref{eiwk2.1} on every open convex subset of ${\cal I}_+$.
\item The set $I_0$ is finite or countable. All limit points of this set belong to ${\cal I}_-$.
\item For any $t_0\in I_0$ Eq. \eqref{eiwk2.2} holds true.
\item The function ${\bar{ z}}(t)$ is a solution of
$$\dot x_k=y_k;\qquad \dot y_k=f_k(t,0,\bar{ z},\mu), \quad k=2,\ldots, n$$
on every open convex subset of ${\cal I}_-$.
\end{enumerate}

Generally speaking, solutions with infinitely many impacts are not unique. For example, one may take the single degree-of-freedom equation $\ddot x=-1$ with delimiter at $x=0$ and $r<1$. The solution $x=0$ is non-unique "backwards"{}.

We identify the vibro-impact system with the pair $({ f},r)$.

Introduce the topology on the set ${\cal X}={\cal X}(J,n,T)$ of vibro-impact systems, corresponding to fixed $J$, $n$ and $T$. This is the minimal topology such that for all $({ f}_0,r_0)\in {\cal X}(J,n,T)$
and every $R>0$ the set
$$\begin{array}{c}\left\{({ f},r)\in {\cal X}:\sup_{(t,{ z},\mu)}\left(|{ f}(t,{ z},\mu)-{ f}_0(t,{ z},\mu)|+\right.\right.\\
+\left|\dfrac{\partial { f}}{\partial { z}}(t,{ z},\mu)-\dfrac{\partial { f}_0}{\partial { z}}(t,{ z},\mu)\right|
+\left|\dfrac{\partial { f}}{\partial \mu}(t,{ z},\mu)-\dfrac{\partial { f}_0}{\partial \mu}(t,{ z},\mu)\right|+\\
\left.\left.|r(\mu)-r_0(\mu)|+|r'(\mu)-r'_0(\mu)|\right)<R\right\}\end{array}$$
is open. The space $\cal X$ with this topology is Hausdorff.

\bigskip

\noindent\textbf{3. The family of periodic solutions}

\bigskip

\noindent Since the solutions of the vibro-impact systems are discontinuous at the impact instants, the classical results on integral continuity are not applicable. Nevertheless, the following statement is true.

\begin{lemma}\label{liwk3.1} Let ${ z}(t)={\mathop {\rm col}\nolimits} (x_1(t),y_1(t),\ldots,x_n(t),y_n(t))$ be the solution of $\eqref{eiwk2.3}$ for $\mu=\mu^*$ and initial conditions ${ z}(t^*_0)={ z}^*_0={\mathop {\rm col}\nolimits} (x^*_{10},y^*_{10},\ldots,x^*_{n0},y^*_{n0})$ such that $x^*_{10}\neq 0$. Suppose that this solution is defined on the segment $[t_-,t_+]\ni t_0$. Assume that there are exactly $N$ zeros $t_-<\tau^1<\ldots<\tau^N<t_+$ of the function $x^*_1(t)$ over the segment $[t_-,t_+]$ and $y_1(\tau^j-0)\neq 0$, $(j=1,\ldots,N)$. Then for any $\varepsilon>0$ there exists a neighborhood $U$ of the point
${\mathop {\rm col}\nolimits} ({ z}^*_0,\mu^*)$ such that for any fixed
$t\in [t_-,t_+]\setminus\bigcup\limits_{k=1}^N (\tau^k-\varepsilon,\tau^k+\varepsilon)$ the mapping ${ z}(t,t_0,{ z}_0,\mu)$ is $C^1$~-- smooth with respect to the variables $(t_0,{ z}_0,\mu)\in U$. Moreover, these solutions have exactly $N$ impact instants
$\tau^j(t_0,{ z}_0,\mu)$ $(j=1,\ldots,N)$ over the segment $[t_-,t_+]$. These instants and corresponding velocities
$$-y_1(\tau^j(t_0,{ z}_0,\mu_1)-0,t_1,{ z}^1,\mu_1)$$
$C^1$~-- smoothly depend on $t_0,{ z}_0$ and $\mu_0$.
\end{lemma}

\textbf{Proof.} Let the number $k$ be such that $t_0\in [\tau^k,\tau^{k+1}]$ (assume, if necessary $\tau^0=t_-$, $\tau^{N+1}=t_+$). The solution ${ z}(t)$ of \eqref{eiwk2.3} is also one of \eqref{eiwk2.1} over any free flight segment. The impact instants $\tau^k$ and $\tau^{k+1}$ as well as the impact points for solutions ${ z}(\tau^k-0,t_0,{ z}_0,\mu)$ and ${ z}(\tau^{k+1}-0,t_0,{ z}_0,\mu)$ smoothly depend on their parameters. Similarly, values $\tau^{k-1}$ and
${ z}(\tau^{k-1}-0,t_0,{ z}_0,\mu)$ smoothly depend on $\tau^k$ and
${ z}(\tau^k-0,t_0,{ z}_0,\mu)$, as well as $\tau^{k+2}$ and ${ z}(\tau^{k+2}-0,t_0,{ z}_0,\mu)$ are $C^1$~-- smooth functions of $\tau^{k+1}$ and
${ z}(\tau^{k+1}-0,t_0,{ z}_0,\mu)$ and so on. $\square$

Suppose that for $\mu\geqslant \mu_0$ the considered system has a family of periodic motions, that pass near the delimiter, and touch it if and only if $\mu=\mu_0$ (Fig.\, 1). More precisely, the following condition is satisfied.

\textbf{Condition 2.} \emph{There exist a segment $\hat J=[\mu_0,\mu_1]\subset J$ and a family of $T$~-- periodic solutions
$${ \varphi}(t,\mu)={\mathop {\rm col}\nolimits}(\varphi_{x1}(t,\mu),\varphi_{y1}(t,\mu),\ldots,\varphi_{xn}(t,\mu),\varphi_{yn}(t,\mu)),\qquad t\in{\mathbb R},\quad \mu\in \hat J$$ of $\eqref{eiwk2.3}$ with the following properties.
\begin{enumerate}
\item For any pair $(t_0,\bar \mu)\in [0,T)\times \hat J$, such that $\varphi_{x1}(t_0,\bar \mu)\neq 0$ the function ${ \varphi}(t,\mu)$ is continuous in a neighborhood of the point $(t_0,\bar \mu)$.
\item For any $\mu\in \hat J$ the component $\varphi_{x1}(t,\mu)$ has $N$ zeros $\tau^1(\mu)<\ldots<\tau^{N}(\mu)$ over the period $[0,T)$. The values $\tau^k(\mu)$ and the velocities
    $$Y^k(\mu)=-\varphi_{y1}(\tau^k(\mu)-0,\mu),\qquad k=1,\ldots, N$$
    continuously depend on $\mu$.
\item For $\mu>\mu_0$ the component $\varphi_{x1}(t,\mu)$ does not have any other zeros. The function $\varphi_{x1}(t,\mu_0)$ has exactly one additional zero $t=0$.
\item Velocities $Y^k(\mu)$ are such that
\begin{equation}\label{eiwk3.1}\begin{array}{c}
\varphi_{y1}(0,\mu_0)=0,
\qquad f_1(0,0,0,\bar{ \varphi}(0,\mu_0),\mu_0)=\phi_0>0,\\
Y^k(\mu)>0, \qquad \forall \mu\in \hat J, \quad k=1,\ldots, N,\\
\min\limits_{\mu\in \hat J}\tau^1(\mu)>0, \qquad \max\limits_{\mu\in \hat J}\tau^N(\mu)<T.
\end{array}\end{equation}
Here $\bar{ \varphi}(t,\mu)={\mathop {\rm col}\nolimits} (\varphi_{x2}(t,\mu),\varphi_{y2}(t,\mu),\ldots,\varphi_{xn}(t,\mu),\varphi_{yn}(t,\mu))$.
\end{enumerate}}

If $\mu$ is the amplitude of free periodic oscillations of the system, the critical value $\mu_0$ corresponds to the distance between the neutral position of the particle and the delimiter. However, in this case, to fulfill Condition 2 we must replace $\mu$ with $-\mu$.

Later on we may suppose without loss of generality that $\bar{ \varphi} (0,\mu_0)=0$, $\mu_0=0$.

\begin{figure}\begin{center}
\includegraphics*[width=1.2in]{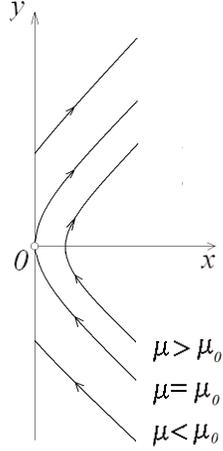}
\end{center}
\caption{The grazing family of periodic solutions.}
\end{figure}

Denote
$$\begin{array}{c}
{ \varphi}_x(t,\mu)={\mathop {\rm col}\nolimits}(\varphi_{x1}(t,\mu),\ldots,\varphi_{xn}(t,\mu)), \\
{ \varphi}_y(t,\mu)={\mathop {\rm col}\nolimits}(\varphi_{y1}(t,\mu),\ldots,\varphi_{yn}(t,\mu)),\\
\bar{{ \varphi}}_x(t,\mu)={\mathop {\rm col}\nolimits}(\varphi_{x2}(t,\mu),\ldots,\varphi_{xn}(t,\mu)), \\
{\bar \varphi}_y(t,\mu)={\mathop {\rm col}\nolimits}(\varphi_{y2}(t,\mu),\ldots,\varphi_{yn}(t,\mu)).\end{array}$$
Consider the shift mapping for system $\eqref{eiwk2.3}$ given by the formula
$${ S}_{\mu}({ z}_0)={ z}(T+ \theta,-\theta,{ z}_0,\mu).$$
Here the value $\theta=\theta(\mu)$ will be specified later (proof of Lemma \ref{liwk7.1}).

For small positive $\mu$ the mapping ${ S}_{\mu}$ is $C^1$~- smooth in a neighborhood of its fixed point
${ z}_{\mu}={ \varphi}(-\theta,\mu)=(x_{1,\mu},y_{1,\mu},{\bar z}_{1,\mu})$.

\bigskip

\noindent\textbf{4. Separatrix}

\bigskip

\noindent Denote $\Gamma_{\mu}=\{{ z}_0\in \Lambda:\exists t_1\in [-T,T]: { z}_1(t_1,0,{ z}_0,\mu)=0\}.$

\begin{lemma}\label{liwk4.1}
There exists a neighborhood $U_0$  of zero such that if the parameter $\mu$ is small enough, the set $\Gamma_{\mu}\bigcap U_0$ is a $2n-1$ dimensional surface, that is the graph of the $C^1$ smooth function $x_1=\gamma_{\mu}(\bar{ x},{ y})$ (Fig.\,2). Moreover,
\begin{equation}\label{eiwk4.1}
\gamma_{\mu}(\bar{ x},{ y})=y^2_1\left(\dfrac1{f_1(0,{ z}_{\mu},\mu)}+\widetilde{\gamma_{\mu}}(\bar{ x},{ y})\right),
\end{equation}
where $\widetilde{\gamma_{\mu}}$ is a $C^1$ smooth function such that $\widetilde{\gamma_{\mu}}(0)=0$.
\end{lemma}

\begin{figure}[h]\begin{center}
\includegraphics*[width=1.2in]{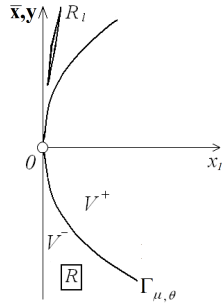}
\end{center}
\caption{The near-grazing behavior of solutions.}
\end{figure}

\textbf{Proof}. Take a point ${ p}_0\in \Gamma_{\mu}$. Let $t_0$ be such that ${ z}_1(t_0,0,{ p}_0)=0$, $s=t-t_0$,
$${ z}_1(t+0,0,{ p}_0)={\mathop {\rm col}\nolimits}(x_1(t),y_1(t)).$$
Let us show that if $t_0$ is close enough to $0$, we may take $s_0\geqslant |t_0|$ so that the function $x_1(t_0+s)$ does not have any zeros on $[-s_0,s_0]$, except $s=0$. Otherwise, there exists a sequence $t^0_k\to 0$ (suppose without loss of generality, that $t^0_k>0$ and the sequence decreases), a sequence
$t^1_k\in [0, t^0_k)$ and one, consisting of solutions:
$${\mathop {\rm col}\nolimits} ({ z}_{k1}(t),\ldots, { z}_{kn}(t))={\mathop {\rm col}\nolimits}(x_{k1}(t),y_{k1}(t),\ldots,x_{kn}(t),y_{kn}(t))$$
of the system $\eqref{eiwk2.3}$, such that ${ z}_{k1}(t^0_k)=0$, $x_{k1}(t^1_k)=0$ (Fig.\, 3).

\begin{figure}\begin{center}
\includegraphics*[width=1.6in]{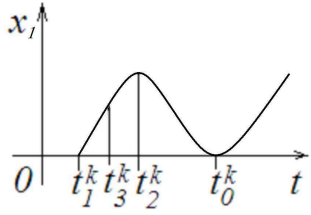}
\end{center}
\caption{Sequences $t^i_k$.}
\end{figure}

Also, there exist time instants $t^2_k\in (t^1_k,t^0_k)$ such that
$\dot x_{k1}(t^2_k)=0$ and instants $t^3_k\in (t^2_k,t^0_k)$ such that $\ddot x_{k1}(t^3_k)=0$. Moreover, $t^3_k\to 0$, $x_{k1}(t^3_k)\to 0$,
$\dot x_{k1}(t^3_k)\to 0$. Without loss of generality, we assume that $\bar{ z}(t_k^3)\to \bar{ z}_0$. Then $\ddot x_1(t^3_k)\to f_1(0,0,\mu)=0$. This contradicts \eqref{eiwk3.1}. Then for all $s\in [-s_0,s_0]$ the function $x_1(t_0+s)$ can be represented as series
\begin{equation}\label{eiwk4.2}
x_1(t_0+s)=X_2s^2+X_3s^3+\ldots
\end{equation}

Differentiating \eqref{eiwk4.2}, we obtain that $\dot x_1(t_0+s)=2X_2s+3X_3s^2+\ldots$. On the other hand,
$X_2=\ddot x_1(t_0+0)/2\to f_1(0,{ z}_{\mu},\mu)$ as $t_0\to 0$. Then
$$\begin{array}{c} x_1(0)=f_1(0,{ z}_{\mu},\mu)(t_0)(1+o(1))/2;\\
y_1(0)=f_1(0,{ z}_{\mu},\mu)(t_0)(1+o(1)).\end{array}$$
Since $y_1=\dot x_1$, \eqref{eiwk4.1} is true. $\square$

Take a small parameter $\varsigma>0$ such that the sets
$$V_{\mu}=\{{ z}\in \Lambda: \|{ z}-{ z}_{\mu}\|\leqslant \varsigma\}\subset U_0,$$
  $$V_{\mu}^-=\{({ z}\in V_{\mu}: x_1<\gamma_{\mu}(\bar{ x},{ y})\}, V_{\mu}^+=\{({ z}\in V: x_1>\gamma_{\mu}(\bar{ x},{ y})\}$$
are correctly defined and nonempty. Here $\|\cdot\|$ is the Euclidean norm.

\bigskip

\noindent\textbf{5. The main result}

\bigskip

\noindent Consider the $(2n)\times(2n)$ matrix
$${ A}=(a_{ij})=({ A}_1,\ldots,{ A}_{2n})=\lim_{\mu\to 0+}{\mathrm D}{ S}_{\mu}({ z}_{\mu}).$$
Suppose the following condition is satisfied.

\textbf{Condition 4.} \emph{The matrix ${ A}$ does not have any eigenvalues on the unit circle. At least one of the eigenvalues of $A$ is inside the unit ball and at least one is outside this ball.}

Note that the matrix $ A$ corresponds to the motion out of a neighborhood of grazing. This may be a motion, described by a linear system without any impacts.

Let $M^s$ be the linear hull of eigenvectors, corresponding to eigenvalues inside the unit ball and $M^u$ be the space, corresponding to eigenvalues outside the unit ball. Let $k=\dim M^s$, $\pi_1$ be the hyperplane of the delimiter, defined by the equality $x_1=0$.

\textbf{Condition 5.} $M^s,M^u\not\subset\pi_1$.

Then the intersections of these spaces with a hyperplane $\pi_1$ are transverse. Denote $\pi_1^{s,u}=\pi_1 \bigcap M^{s,u}$. Let $\alpha_{ij}$ be entries of the matrix $A^{-1}$.
Select a basis ${ e}^s_1, { e}^s_2,\ldots,{ e}^s_k$ in the space $M^s$ and one
${ e}^u_1, { e}^u_2,\ldots,{ e}^u_{2n-k}$
in the space $M^u$ so that ${ e}^s_1 \bot \pi_1^s$, ${ e}^u_1 \bot \pi_1^u$; ${ e}^s_j \in \pi_1^s$, ${ e}^u_j \in \pi_1^u$ for $j>1$.
Denote the components of vectors ${ e}^\sigma_j$ by ${ e}^\sigma_{ij}$, $\sigma\in \{s,u\}$, $i,j=1,\ldots, 2n$. Both the values ${ e}^\sigma_{11}$ are nonzero. Denote
$$R^s={ e}^s_{21}/{ e}^s_{11},\quad R^u={ e}^u_{21}/{ e}^u_{11},\quad R^a=a_{22}/a_{12},\quad R^{\alpha}=\alpha_{22}/\alpha_{12}.$$

\textbf{Condition 6.} \emph{Either
\begin{equation}\label{eiwk5.1}
a_{12}>0, \qquad (R^u-R^a)/(R^a-R^s)>0,
\end{equation}
or
\begin{equation}\label{eiwk5.2}
\alpha_{12}>0,\qquad (R^s-R^\alpha)/(R^\alpha-R^u)>0.
\end{equation}}

From the geometric point of view, Condition 4 means that the point ${ z}_{\mu}$ has nontrivial stable and unstable manifolds ($W^s$ and $W^u$ respectively). Condition 5 means that for small $\mu$ these manifolds intersect transversally the surface $\Gamma_{\mu}$ and bend at the points of intersection. Condition 6 means that the "prolongation"\ of the manifold $W^u$ beyond the intersection with $\Gamma_{\mu}$ intersects transversally with $W^s$ or vice versa (Fig.\, 4).

Later on we suppose that \eqref{eiwk5.1} takes place. Otherwise, if \eqref{eiwk5.2} is true, we consider the mapping ${ S}^{-1}_{\mu}$ instead of ${ S}_{\mu}$. In the proof the matrix $ A$ and the corresponding elements should be replaced with the matrix ${ A}^{-1}$ and corresponding elements and all the references to \eqref{eiwk5.1} should be replaced with ones to \eqref{eiwk5.2}.

\begin{figure}\begin{center}
\includegraphics*[width=3.6cm]{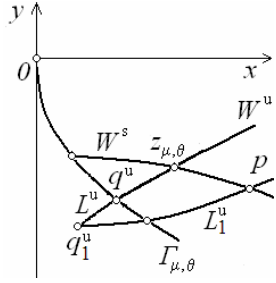}
\end{center}
\caption{Homoclinic points.}
\end{figure}

\begin{theorem}\label{thiwk1} Let Conditions 2 and 4---6 be satisfied. Then there exists a value $\mu_0>0$, such that for all $\mu\in (0,\mu_0)$ the mapping ${ S}_{\mu}$ has a Devaney chaotic invariant set \cite{iwk21}.
\end{theorem}

Namely, there is a hyperbolic transitive infinite invariant set $K_\mu$ of the mapping $S_\mu$ such that periodic points of $S_\mu$ are dense in $K_\mu$.

\bigskip

\noindent\textbf{6. Grazing}

\bigskip

\noindent Now we start proving Theorem \ref{thiwk1}. Fix a small value $\mu>0$ and a solution
$${ z}^*(t)={\mathop {\rm col}\nolimits}(x^*_1(t),y^*_1(t),\ldots,x^*_n(t),y^*_n(t))$$
of the corresponding system, having an impact at $t=t_0$. Suppose that the corresponding normal velocity $Y_{01}=-y_1^*(t_0-0)$ is nonzero but small. Denote $\bar{ Z}_0=\bar{ z}^*(t_0-0)$. Fix a positive value $s_0$ and consider the mapping ${ G}({ q_0})={ z}(t_0+s_0,t_0-s_0,{ q}_0,\mu)$, defined in a neighborhood of the point ${ q}^*_0={ z}^*(t_0-s_0)$. Here we assume that the point ${q}^*_0$ and the parameter $s_0$ are chosen so that there exists a neighborhood $\Omega\ni { q}^*_0$ such that any solution ${ z}(t)={ z}(t,t_0-s_0,{ z}_-,\mu)$ (${ z}_-\in \Omega$) has exactly one impact over the segment $[t_0-s_0,t_0+s_0]$. Denote the corresponding instant by $t_1=t_1({ z}_-)$ and the normal velocity, defined similarly to $Y_{01}$, by $Y_1=Y_1({ z}_-)$. Let $\bar{ Z}=\bar{ z}(t_1)$ be the tangent component of the solution ${ z}(t)$ at the impact instant. Take the values $s_\pm=s_\pm({ z}_-)$ so that $t_0\pm s_0=t_1({ z}_-)\pm s_\pm({ z}_-)$ for all ${ z}_-\in \Omega$. The mapping $ G$ is smooth in the neighborhood of the point $q_0^*$, let us estimate the Jacobian matrix ${\mathrm D}{ G}({ q}_0^*)$. Let
$$\begin{array}{c}
z_-={\mathop {\rm col}} (x_{1,-},y_{1,-},\ldots,x_{n,-},y_{n,-}),\\
x_-={\mathop {\rm col}} (x_{1,-},\ldots,x_{n,-}),\\
y_-={\mathop {\rm col}} (y_{1,-},\ldots,y_{n,-}),\\
{\bar x}_-={\mathop {\rm col}} (x_{2,-},\ldots,x_{n,-}),\\
{\bar y}_-={\mathop {\rm col}} (y_{2,-},\ldots,y_{n,-})
\end{array}$$
Denote
$$\begin{array}{c}
{ z}_+={ z}(t_0+s_0)={ z}(t_0+s_0,t_1+0,0,rY_1,\bar{ Z},\mu),\\
{ x}_+={ x}(t_0+s_0)={ x}(t_0+s_0,t_1+0,0,rY_1,\bar{ Z},\mu),\\
{ y}_+={ y}(t_0+s_0)={ y}(t_0+s_0,t_1+0,0,rY_1,\bar{ Z},\mu),\\
\end{array}$$

Similarly, we define the values $x_{1,+}$, $y_{1,+}$, ${\bar { x}}_{+}$, $\bar { y}_{+}$. Consider the Taylor formula for values
${ z}_\pm$ as functions of $s_\pm$:
\begin{equation}\label{eiwk6.1}\begin{array}{c}
x_{1,-}=Y_1s_-+f_1(t_1,0,-Y_1,{\bar{ Z}},\mu)s_-^2/2+\rho_{x,-}(s_-,t_1,Y_1,{\bar{ Z}},\mu)s_-^2;\\
y_{1,-}=-Y_1-f_1(t_1,0,-Y_1,{\bar{ Z}},\mu)s_-+\rho_{y,-}(s_-,t_1,Y_1,{\bar{ Z}},\mu)s_-;\\
{\bar{ x}}_-={\bar{ x}}(t_1)-{\bar{ y}}(t_1-0)s_-+{\bar f}(t_1,0,-Y_1,{\bar{ Z}},\mu)s_-^2/2+\rho_{{\bar{ x}},-}(s_-,t_1,Y_1,{\bar{ Z}},\mu)s_-^2;\\
{\bar{ y}}_-={\bar{ y}}(t_1-0)-{\bar f}(t_1,0,-Y_1,{\bar{ Z}},\mu)s_-+\rho_{{\bar{ y}},-}(s_-,t_1,Y_1,{\bar{ Z}},\mu)s_-;\\
x_{1,+}=rY_1s_++f_1(t_1,0,rY_1,{\bar{ Z}},\mu)s_+^2/2+\rho_{x,+}(s_+,t_1,Y_1,{\bar{ Z}},\mu)s_+^2;\\
y_{1,+}=rY_1+f_1(t_1,0,rY_1,{\bar{ Z}},\mu)s_++\rho_{y,+}(s_+,t_1,Y_1,{\bar{ Z}},\mu)s_+;\\
{\bar{ x}}_+={\bar{ x}}(t_1)+{\bar{ y}}(t_1-0)s_++{\bar f}(t_1,0,rY_1,{\bar{ Z}},\mu)s_+^2/2+\rho_{{\bar{ x}},+}(s_+,t_1,Y_1,{\bar{ Z}},\mu)s_+^2;\\
{\bar{ y}}_+={\bar{ y}}(t_1-0)+{\bar f}(t_1,0,rY_1,{\bar{ Z}},\mu)s_++\rho_{{\bar{ y}},+}(s_+,t_1,Y_1,{\bar{ Z}},\mu)s_+.\\
\end{array}\end{equation}
Here all functions, denoted by the letter $\rho$, are $C^1$ smooth with respect to all arguments except $s_\pm$ and vanish as $s_\pm=0$. Denote
$$f_{0k+}=f_k(t_0,0,rY_{01},{\bar{ Z}}_0,\mu), \qquad f_{0k-}=f_k(t_0,0,-Y_{01},{\bar{ Z}}_0,\mu).$$
It follows from \eqref{eiwk6.1} that
$$\begin{array}{c}\left. \frac{\displaystyle \partial { z}_+}{\displaystyle \partial (s_+,Y_1,\bar{ Z})}\right|_{s_+=0,\, Y_1=Y_{01},\, \bar{ Z}=\bar{ Z}_0}=
  \begin{pmatrix}
    rY_{01} & 0 & 0\\
    f_{01+} & r & 0\\
    { Q}_+& 0 & { E}_{2n-2}
  \end{pmatrix};  \\[5pt]
  \left. \frac{\displaystyle \partial { z}_-}{\displaystyle \partial (s_-,Y_1,\bar{ Z})}\right|_{s_-=0,\, Y_1=Y_{01},\, \bar{ Z}=\bar{ Z}_0}=
  \begin{pmatrix}
    Y_{01} & 0 & 0\\
   -f_{01-} & -1 & 0\\
   { Q}_- & 0 & { E}_{2n-2}
  \end{pmatrix}.\end{array}$$
Here ${ E}_{2n-2}$ is the unit matrix of the corresponding size,
$${ Q}_+={\mathop {\rm col}\nolimits}(y_2(t_0-0), f_{02+}, \ldots,y_n(t_0-0), f_{0n+}),$$
 ${ Q}_-={\mathop {\rm col}\nolimits}(-y_2(t_0-0), -f_{02+}, \ldots,-y_n(t_0-0), -f_{0n+})$.

Denote
$$f'_k=\left.\dfrac{\partial f_k(t_0,0,y_1,\bar{ Z}_0,\mu)}{\partial y_1}\right|_{y_1=0}.$$
Clearly, $ds_+/ds_-=-1$. Then, similar to the results of the paper \cite{iwk9}, we obtain
\begin{equation}\label{eiwk6.2}{ B}=\lim_{s_\pm \to 0}\frac{\partial{ z}_+}{\partial { z}_-}=\begin{pmatrix}
-r & 0 & 0 \\
b_{21} & -r & 0\\
\bar{{ B}}_1& 0& { E}_{2n-2}
\end{pmatrix}\end{equation}
(Fig.\,2). Here $\bar{{ B}}_1={\mathop {\rm col}\nolimits}(b_{31}, \ldots b_{2n\,1}),$ $$b_{21}=-(f_{01+}+rf_{01-})/Y_{01}=-(r+1)\phi_0(1+O(Y_{01}))/Y_{01},$$ $b_{2j-1\,1}=0$, $b_{2j\,1}=(f_{0j+}-f_{0j-})/Y_{01}=(r+1)f'_k+O(Y_{01})$.
Note that $\det { B}= r^2$.

\bigskip

\noindent\textbf{7. Homoclinic point}

\bigskip

\noindent Clearly, the mapping ${ S}_{\mu}$ is smooth at the points of the set $V_{\mu}$, except ones of the curve $\Gamma_{\mu}$. The matrix ${\mathrm D}{ S}_{\mu}({ z}_0)$ is of the form ${ A}+{ R}_+({ z}_0,\mu)$ where
$$\sup_{{ z}_0\in V^+_{\mu}}{ R}_+({ z}_0,\mu)\to 0$$ as $\mu\to 0+$.
For ${ z}_0\in V^-_{\mu}$ there exist values $t_0({ z}_0)$ and $Y_{01}=Y_{01}({ z}_0)$, continuously depending on ${ z}_0$ and such that
$$x_1(t_0,0,{ z}_0,\mu)=0,\qquad y_1(t_0-0,0,{ z}_0,\mu)=-Y_{01},$$
where $Y_{01}({ z}_0)\to 0$ as ${\mathop {\rm dist}\nolimits} ({ z}_0,\Gamma_{\mu})\to 0$ and
$$\lim_{\mu\to 0+}\sup_{{ z}_0\in V^-_{\mu}}|t_0({ z}_0)|=0.$$

The matrix ${ D}={\mathrm D}{ S}_{\mu}({ z}_0)$ can be represented in the form
${ A}_{\mu}({ z}_0){ B}_{\mu}({ z}_0)$. Here
$${ A}_{\mu}({ z}_0)={ A}+{ R}_-({ z}_0,\mu),\qquad \sup_{{ z}_0\in V^-_{\mu}}{ R}_-({ z}_0,\mu)\to 0$$
as $\mu\to 0+$. The matrix ${ B}_{\mu}({ z}_0)$ can be represented in the form
$${ B}_{\mu}({ z}_0)=({ E}+{ R}_B({ z}_0,\mu)){ B}_{0,\mu}({ z}_0),$$
where ${ R}_B({ z}_0,\mu)=O(1)$ as $\mu\to 0+$ uniformly by ${ z}_0\in V^-_{\mu}$ and the matrix ${ B}_{0,\mu}({ z}_0)$ has the form \eqref{eiwk6.2}. Suppose
${ A}'={ A}_{\mu}({ z}_0)({ E}+{ R}_B({ z}_0,\mu))$. Denote columns of this matrix by ${ A}'_j$ $(j=1,\ldots, 2n)$. Then
$$\begin{array}{l}{ D}({ z}_0)={ A}_{\mu}({ z}_0){ B}_{\mu}({ z}_0)=\\
(-(r+1){ A}'_2\phi_0(1+O(Y_{01}))/Y_{01},\\
-r{ A}'_2(1+O(Y_{01})), { A}'_3+O(Y_{01}), \ldots,{ A}'_{2n}+O(Y_{01})).\end{array}$$

Since $\det { D}({ z}_0)=\det { A}_{\mu}({ z}_0) \det { B}_{\mu}({ z}_0)=(r^2+O(\mu))\Delta_0$ and $a_{12}\neq 0$, one of the eigenvalues of the matrix ${ D}({ z}_0)$ is equal to
$$\lambda_+({ z}_0)=-(r+1)a_{12}\phi_0(1+O(Y_{01}))/Y_{01}.$$
The corresponding eigenvector is of the form ${ u}_+={ A}_2+O(\mu)$. The linear space, corresponding to other eigenvalues, tends to the hyperplane $\pi_1=\{(0,y_1,\bar z)\}$ as $\mu\to 0$. Since $a_{12}\neq 0$, the vector ${ u}_+({ z}_0)$ does not belong to $\pi_1$. The vector ${ A}_2$ (as well as the vector ${ A}'_2$) is out of the space $M^s$ and the space $M^u$.

The matrix ${ D}^{-1}({ z}_0)$ satisfies the following asymptotic estimates:
$${ D}^{-1}=\dfrac1{r^2}\begin{pmatrix}
-r{ A}^-_1+O(\mu)\\
(r+1){ A}^-_1(1+O(\mu))/Y_{01}\\
{ A}^-_3+O(\mu)\\
\dots\\
{ A}^-_{2n}+O(\mu)
\end{pmatrix},$$
where ${ A}^-_j$ are the strings of the matrix ${ A}^{-1}$. Hence there is an eigenvalue
$$\lambda_-^{-1}=(r+1)\phi_0 \alpha_{12}(1+O(\mu))/(r^2Y_{01})$$ of the matrix ${ D}^{-1}$. The corresponding eigenvector is $u_-={ e}_{2}+O(\mu)$. Here $${ e}_2={\mathop{\rm col}\nolimits} (0,1,0,\ldots,0).$$

It follows from the Perron theorem that in a small neighborhood of the point ${ z}_{\mu}$ there exist local stable and unstable manifolds  $W^s_{loc}$ and $W^u_{loc}$ respectively. Both of them are smooth in a neighborhood of the point ${ z}_{\mu}$. Clearly,
$M^{s,u}=T_{{ z}_{\mu}} W^{s,u}_{loc}$. Extend the manifolds $W^s_{loc}$ and $W^u_{loc}$ up to invariant sets of the mapping ${ S}_{\mu}$. Denote these sets by $W^s$ and $W^u$, respectively. Generally speaking, both of these sets consist of a countable number of connected components. All these components are piecewise smooth manifolds.

\begin{lemma}\label{liwk7.1} The intersection of the sets $W^s$ and $W^u$ contains a point ${ p}\neq { z}_{\mu}$ (Fig.\, 4). The connected components of intersections of the sets $W^s$ and $W^u$ with a neighborhood of the point $ p$, containing this point, are smooth manifolds and their intersection is transverse.
\end{lemma}

\textbf{Proof.} If the parameter $\mu$ is small enough, the manifold $W^u$ transversally intersects with the surface $\Gamma_{\mu}\bigcap V_{\mu}$. Denote the $k-1$~-- dimensional smooth manifold, obtained in the intersection, by $q^u$.

Let $q^u_1={ S}_{\mu}(q^u)$. The neighborhood $V_{\mu}$ of the point ${ z}_{\mu}$ may be chosen so that both the manifolds $q^u$ and $q^u_1$ intersect with  $V_{\mu}$. Let $L^u$ be the connected component of the intersection of $(W^u\bigcap V_{\mu}) \setminus (q^u\bigcup q_1^u)$ whose boundary intersects with $q^u$ and $q^u_1$. We select the parameter $\theta(\mu)$ in the definition of $S_\mu$ so small that $q^u_1$ and $L^u$ are correctly defined.

The neighborhood $V_{\mu}$ may be chosen so that
${\mathop {\rm diam}\nolimits}\, V_{\mu}\to 0$ as $\mu\to 0+$. Let $L^u_1={ S}_{\mu}(L^u)$. For any ${ z}\in L^u$ the tangent space $M^u({ z})=T_z W^u$ is the linear hull of unit orthogonal vectors ${ v}^u_1({ z}),{ v}^u_2({ z}), \ldots, { v}^u_{2n-k}({ z})$, taken so that for all $k=1,\ldots,n_u$
$$\lim_{\mu\to 0+}\max_{z\in L^u}\|{ v}^u_k({ z})-{ e}^u_k\|=0.$$
Recall that vectors $e^u_k$ have been defined immediately after Condition 5.

The surface $W^u$ is not smooth in a neighborhood of the manifold $q^u_1$. For the points ${ z}\in L_1^u$ the tangent space $M^u({ z})=T_{ z} W^u$ is a linear hull of unit vectors $\widetilde{{ v}}^u_1({ z}),{ v}^u_2({ z}), \ldots, { v}^u_{n_u}({ z})$ such that
$$\lim_{\mu\to 0+}\max_{z\in L^u_1}\|{ v}^u_k({ z})-{ e}^u_k\|=0$$
for all $k>1$ and
$$\lim_{\mu\to 0+}\max_{z\in L^u_1}\|\widetilde{v}^u_1({ z})-{ A}_2/\|{ A}_2\|\|=0.$$
Here we use the fact that $L_1^u$ contains the manifold $q_1^u$ whose inclination with respect to $\pi_1$ is arbitrarily small if $\mu$ is small.

It follows from Condition 6 that for all ${ z}\in L^u_1$ vectors
$$\widetilde{ v}^u_1({ z}),{ v}^u_2({ z}), \ldots, { v}^u_{n_u}({ z})$$
are linearly independent. Vectors ${ A}_2$ and ${ e}^u_1$ belong to different half-spaces, separated by the hyperplane $\pi_1$. Consequently, for all ${ z}_0\in L^u$, ${ z}_1\in L^u_1$ and ${ z}_2\in W^u\bigcap V_{\mu}$ vectors ${ v}^u_1({ z}_0)$ and $\widetilde{v}^u_1({ z}_1)$ belong to different half-spaces, separated by the linear hull of vectors
$${ v}^u_2({ z}_2), \ldots, { v}^u_{n_u}({ z}_2),{ v}^s_1({ z}_2), \ldots, { v}^s_{n_s}({ z}_2).$$

Let $p_\mu$ be the nearest to ${ z}_{\mu}$ point of the set $\pi_1 \bigcap \widetilde{M}^u$. The second coordinate of ${ p}^u$ equals to
\begin{equation}\label{eiwk7.1}
y_0=y_{1,\mu}-R^u x_{1,\mu}.
\end{equation}
Recall that $x_{1,\mu}$ and $y_{1,\mu}$ are the first two coordinates of the point ${ z}_{\mu}$.

Condition \eqref{eiwk5.1} is equivalent to the existence of a solution $(X,Y)=(x_p,y_p)$ of the system
$$\dfrac{X-x_{1,\mu}}{Y-y_{1,\mu}}=R^s;\qquad \dfrac{X}{Y-y_0}=R^a$$
such that $x_p>x_{1,\mu}$. Here the value $y_0$ may be found by \eqref{eiwk7.1}. Hence, the affine space, containing the point $p_0$ and parallel to the linear hull of the vectors ${ A}_2$, ${ e}^u_2$, $\ldots$, ${ e}^u_k$ intersects the space, tangent to the manifold $W^u$ at the point ${ z}_{\mu}$. The first coordinate of the intersection point is greater than $x_{1,\mu}$, that provides a transverse intersection of $L^u_1$ and $W^s$, if $\mu$ is small enough. $\square$

The parameter $\mu$ is chosen so that the point of this intersection, nearest to ${ z}_{\mu}$ does not belong to the hyperplane $\pi_1$.

The Smale-Birkhof theorem \cite{iwk21} on the existence of a chaotic invariant set in a neighborhood of a homoclinic point is not applicable in the considered case since the mapping ${ S}_{\mu}$ is discontinuous. So we have to find a chaotic set of the mapping ${ S}_{\mu}$ "manually".

\bigskip

\noindent\textbf{8. Symbolic dynamics}

\bigskip

\noindent Consider new smooth coordinates ${\mathop {\rm col}\nolimits}(\zeta_1^s,\ldots, \zeta_{n_s}^s,\zeta_1^u,\ldots, \zeta_{n_u}^u)$ in the neighborhood of the point ${ z}_{\mu}$ such that the following conditions are satisfied.
\begin{enumerate}
  \item The point ${ z}_{\mu}$ corresponds to zero in the new coordinate system.
  \item The Euclidean norm of any column of the matrix $$\frac{\displaystyle \partial \, {\mathop {\rm col}\nolimits}(\zeta_1^s,\ldots, \zeta_{n_s}^s,\zeta_1^u,\ldots, \zeta_{n_u}^u)}{\displaystyle \partial { z}}({ z}_{\mu})$$ equals to 1.
  \item In a small neighborhood of zero the stable and the unstable manifolds are given by the conditions $\zeta^u_j=0$ and $\zeta^s_j=0$ $(j=1,\ldots,n)$ respectively.
  \item The direction of the tangent line to the axis, $O\zeta^s_1$, coincides with one of the vector ${ u}_+$, and one, corresponding to the axis $O\zeta^u_1$ coincides with the vector ${ u}_-$.
\end{enumerate}
Consider the neighborhood $Q_0$ of the point ${ z}_{\mu}$, defined by conditions
$$|\zeta^s_j|\leqslant \varepsilon_s, \quad |\zeta^u_j|\leqslant \varepsilon_u,\qquad j=1,\ldots,n$$
(Fig.\,5). For any $m\in {\mathbb Z}$ define
$Q_m={ S}_{\mu}^m(Q_0)$. Denote parts of the boundary of the set $Q_0$, that correspond to $\zeta^s_1=\pm\varepsilon_s$ and $\zeta^u_1=\pm\varepsilon_u$ by $\partial^\pm_x$ and $\partial^\pm_y$ respectively. Here the positive values $\varepsilon_s$ and $\varepsilon_u$ are chosen so that $\varepsilon_s\geqslant 2\varepsilon_u$ and there exist natural numbers $m_+$ and $m_-$ such that ${ p}\in Q_{m^+}\bigcap Q_{-m^-}$ (Fig.\,6). We may take $\mu$ and, respectively, eigenvalues $\lambda^+({ z}_{\mu})$ so that $Q_j\subset V^-_{\mu}$ for any $-m^-<j<m^+$.

\begin{figure}\begin{center}
\includegraphics*[width=1.83in]{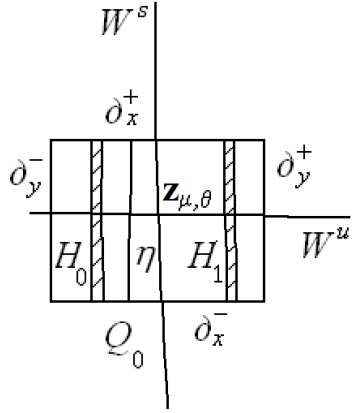}
\end{center}
\caption{Domains $H_0$ and $H_1$.}
\end{figure}

Due to Lemma \ref{liwk7.1} the set $Q_{-m^-}\bigcap Q_{m^+}$ contains at least two connected components. One of them denoted by $\widetilde{H}_0$ contains the point ${z}_{\mu}$, another one denoted by $\widetilde{H}_1$ contains the point ${ p}$. Let $m=m^++m^-$, $H_0={ S}_{\mu}^{m^-}(\widetilde{H}_0)$,
$H_1={ S}_{\mu}^{m^-}(\widetilde{H}_1)$, $H=H_0\bigcup H_1$ (Fig.\,5).


Let us show that the set $$K=\bigcap\limits_{k\in {\mathbb Z}}{ S}_{\mu}^{mk}(H)$$ is chaotic. Evidently, the set $K$ is invariant with respect to the mapping ${ S}_{\mu}$, compact and nonempty, since it contains the point ${ z}_{\mu}$. Moreover, neither the inverse images ${ S}_{\mu}^{-mk}\pi_1$ ($k\in{\mathbb Z}$) of the hyperplane $\pi_1$ nor ones of the surface $\Gamma_{\mu}$ intersect with $K$. Consequently, for any integer $k$ there is a neighborhood $U_{\mu}$ of the set $K$ such that the mapping ${ S}_{\mu}^{m}|_{U_{\mu}}$ is $C^1$ smooth.

\begin{lemma}\label{liwk8.1}
For any $k\in {\mathbb N}$, any set ${ i}=(i_0, \ldots, i_k)$ such that $i_j\in \{0,1\}$ for any $j=0,\ldots, k$, the set $H_i=H_{i_0}\bigcap { S}_{\mu}^{-m}(H_{i_1})\bigcap\ldots\bigcap { S}_{\mu}^{-mk}(H_{i_k})$ is not empty.
\end{lemma}

\textbf{Proof.} Consider an arbitrary arc $\eta$, joining disks $\partial^+_x$ and $\partial^-_x$ and defined as the graph of a smooth function
$$(\zeta^s_1,\ldots,\zeta^s_n,\zeta^u_2,\ldots,\zeta^u_n)={ h}(\zeta^u_1),$$
such that $\max \|{ h}'(\zeta^u_1)\|\leqslant 1$. Let us call such arcs admissible. Similarly to the Palis lemma \cite[chapter 2, Lemma 7.1]{iwk21}, one may show that for small values of $\mu$, $\varepsilon_s$ and $\varepsilon_u$ there exists an embedding of the curve ${ S}_{\mu}^{-m}(\eta)$ to the manifold $W^s$, arbitrarily $C^1$-close to the identical embedding of ${ S}_{\mu}^{-m}(\eta)$ to ${\mathbb R}^{2n}$. Particularly, this means that the set ${ S}_{\mu}^{-m}(\eta)$ contains two admissible arcs $\eta_0$ and $\eta_1$. Fix an index $ i$. It follows from what is proved, that the inverse image ${ S}_{\mu}^{-m}(\eta)$ of any admissible arc $\eta$, contains an admissible arc
$\eta_{i_0}\subset H_{i_0}$. Applying the same procedure to the curve $\eta_{i_0}$, we obtain an arc
$\eta_{i_0i_1}\subset { S}_{\mu}^{-m}(\eta_{i_0})\bigcap H_{i_1}$. Finally, we get a curve
$\eta_{i_0i_1\ldots i_m}\subset { S}_{\mu}^{-mN}(H_{i_0})\bigcap { S}_{\mu}^{-mN+m}(H_{i_1})\bigcap \ldots \bigcap H_{i_m}$.
Then the statement of the lemma follows from the inclusion ${ S}_{\mu}^m(\eta_{i_0i_1\ldots i_m})\subset H_i$.

Hence, for any point ${ z}\in K$ there is a unique sequence $${ i}({ z})=\{\ldots,i_{-2},i_{-1},i_0,i_1,i_2,\ldots\}$$
such that ${ S}_{\mu}^{mk}({ z})\in H_{i_k}$ for any $k\in{\mathbb Z}$. Due to Lemma \ref{liwk8.1} for any sequence ${ i}$ one may find the corresponding point $ z$. Since the diffeomorphism ${ S}_{\mu}^m$ is hyperbolic in a neighborhood of the set $K$ the point ${ z}\in K$ is uniquely defined by the sequence ${ i}({ z})$. Therefore the set $K$ is of the power continuum. The unit shift of the index to the left corresponds to the mapping ${ S}_{\mu}^m$. The presence of this symbolic dynamics proves the theorem. $\square$

\bigskip

\noindent\textbf{9. Example}

\bigskip

\noindent Consider the equation $\ddot x+p\dot x-qx=-a-b\sin(\omega t)$ equivalent to the system
\begin{equation}\label{eiwk9.1}
\left\{\begin{array}{l}
\dot x=y;\\
\dot y=qx-py-a-b\sin(\omega t).
\end{array}\right.
\end{equation}
We suppose that all the parameters of this system are real, and $a$, $q$ and $\omega$ are positive. Consider the vibro-impact system, defined by Eq. \eqref{eiwk9.1} and the one-dimensional impact condition, where $r=1$. The first component $x(t)$ of the general solution of \eqref{eiwk9.1} is of the form
$$x=C_+\exp(\lambda_+ t)+C_-\exp (\lambda_- t)+\varphi(t).$$
Here
$$\lambda_+=\dfrac{-p+\sqrt{p^2+4q}}2>0,\qquad \lambda_-=\dfrac{-p-\sqrt{p^2+4q}}2<0,$$
$$\varphi(t)=\dfrac{a}{q}+A\sin(\omega t)+B\cos(\omega t),$$
where
$$A=\dfrac{b(\omega^2+q)}{(\omega^2+q)^2+p^2\omega^2}, \qquad B=\dfrac{bp\omega}{(\omega^2+q)^2+p^2\omega^2}.$$
Note that $\varphi(t)$ corresponds to the unique periodic solution of \eqref{eiwk9.1}. The function $\varphi(t)$ is positive for all $t$ if and only if
$A^2+B^2<a^2/q^2$ or, equivalently, if $b^2q^2<a^2((\omega^2+q)^2+p^2\omega^2)$. The fundamental matrix of the corresponding homogeneous system, that turns to the unit matrix for $t=0$ is
$\Phi(t)=$
$$\dfrac1{\sqrt{p^2+4q}}\begin{pmatrix}
\lambda_+\exp(\lambda_- t)-\lambda_-\exp(\lambda_+ t) & \exp(\lambda_+ t)-\exp(\lambda_- t) \\
-q(\exp(\lambda_+ t)-\exp(\lambda_- t)) & \lambda_+\exp(\lambda_+ t)-\lambda_-\exp(\lambda_- t).
\end{pmatrix}$$

The matrix $A$ is equal to $\Phi(2\pi/\omega)$. The spaces $M^s$ and $M^u$ are linear hulls of vectors
$\mbox{\rm col}(1,\lambda_-)$ and $\mbox{\rm col}(1,\lambda_+)$ respectively and
$$R^s=\lambda_-, \quad R^u=\lambda_+, \quad R^a=\dfrac{\lambda_+\exp(\lambda_+ t)-\lambda_-\exp(\lambda_- t)}{\exp(\lambda_+ t)-\exp(\lambda_- t)}.$$

Conditions 4 and 5 hold true evidently. Clearly, $a_{12}>0$ for all values of parameters of \eqref{eiwk9.1}. Moreover, for all $T$ we have $R^a>R^u>R^s$, that implies \eqref{eiwk5.1}. So, the conditions of Theorem \ref{thiwk1} are satisfied.

\bigskip

\noindent\textbf{10. "Soft"\ model and structural stability}

\bigskip

In this section we do not need the right hand side ${ f}$ of Eq. \eqref{eiwk2.1} to depend on the parameter $\mu$, so we omit this parameter in our notations.

Suppose that Eq. \eqref{eiwk2.1} is defined for ${ z}\in {\Lambda}$. Fix a parameter value $r\in (0,1]$ and denote Eq. \eqref{eiwk2.3}, corresponding to this fixed value of $r$ by $\eqref{eiwk2.3}_r$. Note, that any Cauchy problem for the system $\eqref{eiwk2.3}_r$ with initial conditions
${ z}(t_0)={ z}_0$, $t_0\in {\mathbb R}$,
$${ z}_0=\mbox{\rm col}\,(x_{10},y_{10},{\bar{ z}_0})\in \Lambda$$
has a solution. It is locally unique if $x_{10}\neq 0$. In this case there is a value $\varepsilon>0$, such that on time intervals $(t_0-\varepsilon,t_0)$ and $(t_0,t_0+\varepsilon)$ the considered solution does not have any impacts and, consequently, coincides with a solution of Eq. \eqref{eiwk2.1}.

Consider a function $\chi_-(s)=1-H(s)$ where $H(s)$ is the standard Heaviside step function, i.e. $\chi_-(s)=1$ if $s<0$ and $\chi_-(s)=0$ if $s\ge 0$. For a fixed value $r\in (0,1]$ define $\alpha=-(\log r)/\pi$. We define the function $f(t,{ z})$ on the set ${\mathbb R}^{2n+1}\setminus \Lambda$ so that this function obtained is piecewise smooth on the whole space ${\mathbb R}^{2n+1}$. Let $\nu>0$ be a big parameter. Define
$$h(\nu,{ z})=(-2\alpha \nu y_{1}-(1+\alpha^2)\nu^2x_1)\chi_-(x_1){ e}_1.
$$
Here ${ e}_1=\mbox{\rm col}\, (1,0,\ldots,0)$. The function $h$ is piecewise continuous.

In this section we study the following problems.
\begin{enumerate}
  \item When do the invariant sets of the shift mapping for Eq. $\eqref{eiwk2.3}_r$ persist provided the coefficients of Eq. \eqref{eiwk2.1} and/or the restitution coefficient $r$ are slightly changed?
  \item When these invariant sets persist provided one replaces the impact condition with the perturbation $h(\nu,{ z})$ corresponding to a big value of $\nu$?
\end{enumerate}

Let $g(t,z)={\mathop {\rm col}} (g_1(t,z),\ldots,g_n(t,z))$ be a $C^1$ smooth mapping, $T$~- periodic with respect to $t$. Suppose that it is small enough in the $C^0$~-- norm together with $\partial g/\partial z$. Consider the following systems
\begin{equation}\label{eiwk10.1}
\dot x_k=y_k, \qquad \dot y_k=f_k(t,{ z})+g_k(t,{ z})\qquad \mbox{and}
\end{equation}
\begin{equation}\label{eiwk10.2}
\dot x_k=y_k, \qquad \dot y_k=f_k(t,{ z})+g_k(t,{ z})+h(\nu, {z}).
\end{equation}

We denote the dynamical system, which consists of Eq. \eqref{eiwk10.1} and Condition 1 where $r$ is fixed by $\eqref{eiwk10.1}_{g,r}$. Consider the solution
$z_r(t,t_0,z_0)$ of the system $\eqref{eiwk2.3}_r$ with initial conditions $z(t_0)=z_0$ and functions $z_{g,\nu,r}(t,t_0,z_0)$ and $z_{g,r}(t,t_0,z_0)$, which are solutions of Eq. \eqref{eiwk10.2} and $\eqref{eiwk10.1}_{g,r}$ with the same initial conditions. Introduce shift mappings for Eq. $\eqref{eiwk2.3}_r$, $\eqref{eiwk10.1}_{g,r}$ and \eqref{eiwk10.2}, by following formulae $G_r(z_0)=z_r(T,0,z_0)$, $G_{g,r}(z_0)=z_{g,r}(T,0,z_0)$, $G_{g,\nu,r}(z_0)=z_{g,\nu,r}(T,0,z_0)$.

\begin{theorem}\label{thiwk2}
Let $r_0\in (0,1]$. Suppose that the mapping $G=G_{r_0}$ has a hyperbolic invariant set $K\subset (0,+\infty)\times {\mathbb R}$ such that $z_{1,r_0}(t,0,z_0)\neq 0$ for all $z_0\in K$, $t\in [0,T]$. Let $U$ be a neighborhood of $K$ such that $\overline{U}$ and the corresponding image $G(\overline{U})$ do not intersect with the axis $Oy$. Let
$$W=\{t,z:t\in [0,T], z=z(t,0,z_0), z_0\in \overline{U}\}.$$
The following statements are true.
\begin{enumerate}
  \item For every $\varepsilon >0$ there exists $\delta>0$ such that if $r\in (r_0-\delta,r_0+\delta)\bigcap (0,1]$,
\begin{equation}\label{eiwk10.3}
\max\limits_{(t,z)\in W} |g(t,z)|<\delta, \quad \max\limits_{(t,z)\in W} \left|\dfrac{\partial g}{\partial z}(t,z)\right|<\delta,
\end{equation}
the mapping $G_{g,r}$ is well-defined in a neighborhood $U_0$ of $K$. There exists a homeomorphism
$h_{g,r}:K\to K_{g,r}\subset U_0$ such that $\max\limits_{x\in K}|h_{g,r}(x)-x|<\varepsilon$, and the set $K_{g,r}$ is hyperbolic invariant for the mapping
$G_{g,r}$. Moreover, for any $x\in K$
\begin{equation}\label{eiwk10.4}
h_{g,r}(G(x))=G_{g,r}(h_{g,r}(x)).
\end{equation}
\item Let $r_0\in (0,1]$. For any $\varepsilon >0$ there exist $\nu_0>0$, $\delta>0$ such that if $\nu>\nu_0$ and conditions \eqref{eiwk10.3} are
satisfied, there exists a homeomorphism $\eta_{g,\nu,r_0}:K\to K_{g,\nu,r_0}\subset {\mathbb R}^2$ such that
$\max\limits_{x\in K}|\eta_{g,\nu,r_0}(x)-x|<\varepsilon$ and $K_{g,\nu,r_0}$ is a hyperbolic invariant set of $G_{g,\nu,r_0}$. Moreover,
$\eta_{g,\nu,r_0}(G(x))=G_{g,\nu,r_0}(\eta_{g,\nu,r_0}(x))$ for all $x\in K$.
\end{enumerate}
\end{theorem}

\textbf{Proof.} Let us check item (1). The neighborhood $U$  can be chosen so that any solution of Eq. $\eqref{eiwk10.1}_{r_0}$, which starts at the instant $t=0$ in the domain $\overline{U}$, has at most a fixed value $M$ of impacts. Let it be not true. Then there exists a sequence $z_k\in
\overline{U}$, $k\in {\mathbb N}$ of initial conditions, corresponding to solutions, which have at least $k$ impacts over $[0,T]$. Without loss of generality, one may assume, that $z_k\to z_0\in \overline{U}$. If such a point exists for any choice of the neighborhood $U$, we may say that $z_0\in K$. The solution $z(t,0,z_0)$ has infinitely many impacts on $[0,T]$, consequently, instants of these impacts have a limit point $t^*\in [0,T]$. Then $z(t^*,0,z_0)=0$. This contradicts to assumptions of the theorem.

The set $\overline{U}$ can be represented in the form $\overline{U}=Q_0\bigcup \ldots \bigcup Q_M,$ where every set $Q_k$ is compact or empty and consists of initial conditions corresponding to solutions of Eq. $\eqref{eiwk10.1}_{r_0}$, which have exactly $m$ impacts on the segment $[0,T]$, such that corresponding values of the normal velocity $y_1$ are always nonzero. Then the impact instants smoothly depend on the initial conditions. On every set $Q_k$ except $Q_0$ the mapping $G$ is of the form $G=I_0\circ G^1 \circ \ldots \circ G^{k-1} \circ I_1$. Here $I_0$ maps the initial conditions $z_0\in \Lambda$ to the first impact instant $t^1$, the point $\bar x(t^1)\in \pi_1$ and the velocity $y(t^1)$, corresponding to the first impact. The mappings $G^j$, $j=1,\ldots,k-1$, transfer the instant, position and velocity of the impact number $j$ into correspondence to the same parameters of the $j+1$-th impact. The mapping $I_1$ transfers the instant, point and velocity of the last impact to the phase coordinates of the corresponding solution at the instant $t=T$. Clearly, mappings $I_0$, $I_1$ and $G^j$ ($j=1,\ldots,k-1$) are $C^1$ smooth on the set $Q_l$. Hence there is a $\sigma>0$ such that the mapping $G$ is $C^1$ smooth on the set $\overline U\times [r_0-\sigma,r_0+\sigma]$.

The following auxiliary statement is analogous to the theorem on persistence of hyperbolic invariant sets of diffeomorphisms \cite{iwk23}. The only difference is that we consider an embedding of a domain instead of a diffeomorphism of a manifold. The proof, given at \cite{iwk23} is still valid for the considered case.

\begin{lemma}\label{liwk10.1} Let $d$ be a natural number, $U\subset {\mathbb R}^d$ be a domain. Let the $C^1$ smooth embedding $\phi:U \to {\mathbb R}^d$ possess a compact hyperbolic invariant set $K\subset U$. Then for any $\varepsilon>0$ there exists a $\delta>0$ such that if the mapping $\psi:U\to {\mathbb R}^d$ is such that $\|\psi-\phi\|_{C^1(U\to {\mathbb R}^d )}<\delta$, then there is an embedding $h:K\to {\mathbb R}^d$ satisfies the inequality
$\|h-{\mathop{\rm id}\nolimits}\|_{C^0(K\to {\mathbb R}^d)}<\varepsilon$ and
\begin{equation}\label{eiwk10.5}
h(\phi(x))=\psi(h(x))
\end{equation}
for all $x\in K$. Particularly, the set $K_1=h(K)$ is hyperbolic invariant for the mapping $\psi$ if $\varepsilon$ is small enough.
\end{lemma}

Applying Lemma \ref{liwk10.1} to the mapping $G_{r_0}$ and its small perturbations $G_{g,r}$, we obtain that there exist $\delta>0$ and a neighborhood $V\subset U$ of the set $K$, such that for any $r\in (r_0-\delta,r_0+\delta)\bigcap (0,1]$ and any perturbation $g$, satisfying conditions \eqref{eiwk10.3}, there exists a homeomorphism $h_{g,r}$, topologically conjugating $G$ and $G_{g,r}$ in the sense of \eqref{eiwk10.4}. Moreover
$h_{g,r}\rightrightarrows {\mathop{\rm id}\nolimits}$ as $\delta\to 0$ and $r\to r_0$.

So, the mapping $G_{g,r}$ is a diffeomorphism and the set $K_{g,r}=h_{g,r}(K)$ is hyperbolic invariant. This proves the first part of the theorem. Let us start to prove the second one.

\begin{lemma}\label{liwk10.2}  For all $y_0^+>y_0^->0$, any $t_0\in {\mathbb R}$ and all functions $f$ and $g$, satisfying conditions of Theorem \ref{thiwk2}, there exists a value $\nu_1>0,$ such that if $\nu>\nu_1$, the solution of Eq. \eqref{eiwk10.2} with initial conditions
$$x(t_0)=x_0\in \pi_1, \qquad y(t_0-0)=y_0=\mbox{\rm col}\,(y_{1,0},{\bar y}_0): y_{1,0}\in(-y_0^+,-y_0^-)$$
intersects $\pi_1$ at $t_1\in (t_0,t_0+2\pi/\nu)$. Moreover there exists a $\Delta>0$, which does not depend on $y_0$, such that for any
$\nu>\nu_1$
\begin{equation}\label{eiwk10.6}
|y_{1}(t_1)+r_0y_{N,0}|<\Delta/\nu.
\end{equation}
\end{lemma}

\textbf{Proof.} The transformation of the independent variable $s=\nu t$ reduces \eqref{eiwk10.2} for $x\notin \bar{\Lambda}$ to the form
\begin{equation}\label{eiwk10.7}\begin{array}{c}
x_{1}''+2\alpha x_{1}'+(1+\alpha^2)x_{1}-(f_{1}(s/\nu,x,\nu x')+
g_{1}(s/\nu,x,\nu x')))/{\nu^2}=0,\\
\bar x''-(\bar f(s/\nu,x,\nu x')+ \bar g(s/\nu,x,\nu x'))/{\nu^2}=0.
\end{array}\end{equation}
Here we denote by prime the derivative $d/ds$. The initial conditions $(t_0,x_0,y_0)$ are transformed to the following ones:
\begin{equation}\label{eiwk10.8}\begin{array}{c}
s_0=\nu t_0,\quad x_{1}(s_0)=0,\quad
\bar x(s_0)=x_0,\\
x_{1}'(s_0)=y_{1,1}/\nu \in(-y_0^+/\nu,-y_0^-/\nu), \quad
\bar x'(s_0)=\bar y_{0}/\nu.
\end{array}\end{equation}
If $\nu$ is big enough, Eq. \eqref{eiwk10.7} on compact sets is a small perturbation of the system
\begin{equation}\label{eiwk10.9}
 x_{1}''+2\alpha x_{1}'+(1+\alpha^2)x_{1}=0, \qquad \bar x''=0.\end{equation}

The first component of the general solution $\phi$ of Eq. \eqref{eiwk10.9} is of the form
$$x_{1}=A\exp(-\alpha s)\sin (s+\varphi).$$
The distance between neighbor zeroes of this function is always equal to $\pi$, consequently, for $\nu$ big enough and for any solution $\Psi\neq 0$ of Eq. \eqref{eiwk10.7} with initial conditions \eqref{eiwk10.8} the distance between neighbor zeros of the first component of $\Psi$ varies from $\pi/2$ to $2\pi$. For Eq. \eqref{eiwk10.2}, we obtain that if the solution with initial conditions on $\overline{U}$ intersects with $\pi_1$ at the instant $t_0$ and the corresponding normal component of velocity is negative. The solution intersects $\pi_1$ again at the instant $t_1\in (t_0,t_0 +2\pi/\nu)$. The ratio of normal components of derivatives of solutions of Eq. \eqref{eiwk10.9} corresponding to successive impact instants is $-r_0$. Since the solutions between impacts continuously depend on initial conditions and parameters, there exists $\nu_2>0$, such that if $\nu>\nu_2$, for any solution $x(t)$ of
Eq. \eqref{eiwk10.7} with initial conditions \eqref{eiwk10.8} there is a $\Delta>0$, such that $|x_{1}'(\nu t_1)+r_0x_{1}'(s_0)|\leqslant \Delta/\nu^2.$ This proves the validity of \eqref{eiwk10.6}. $\square$

Let the number $m\in \{0,\ldots,M\}$ be such that $z_0\in Q_m$. If $m=0$, the function $z_{g,\nu,r}(t,0,z_0)$ is a solution of Eq. \eqref{eiwk10.2}. Else, for $\nu$ big enough and a $\delta>0$ such that $|z_0|<\delta$, the normal component $x_{1,g,\nu,r}(t,0,z_0)$ of the solution $z_{g,\nu,r}(t,0,z_0)$ has exactly $2m$ zeroes on the segment $[0,T]$. Denote them by $t^1<\theta^1<t^2<\theta^2<\ldots <t^k<\theta^k$. The normal component $y_{1,g,\nu,r}(t,0,z_0)$ of the velocity of the solution $z_{g,\nu,r}(t,0,z_0)$ is positive at the instants $t^j$ and negative at the instants $\theta^j$. The
values $\nu$ and $\delta$ may be chosen independently with respect to $z_0$ and $m$. The mapping $G_{g,\nu,r}$ is of the form $G_{g,\nu,r}=\widehat{I}_0\circ H^1\circ {H^1}' \circ H^2\circ {H^2}'\ldots \circ H^{k-1}\circ {H^{k-1}}'\circ \widehat{I}_1$. The mapping $\widehat{I}_{0}$ transfers initial conditions $z_0\in \Lambda$ at the time instant $t=0$ to the triple $(t^1,X^1,Y^1)$. Here $t^1$ is the first zero of the normal component $x_{1,g,\nu,r}(t,0,z_0)$, the tangential component $X^{1}=\bar x_{g,\nu,r}(t^1,0,z_0)$ and the value $Y^1=y_{g,\nu,r}(t^1,0,z_0)$. Mappings $H^j$, $j=1,\ldots,k-1$ transfer triples $(t^j,X^j,Y^j)$ of time instants $t^j$, impact points $X^j\in \pi_1$ and velocities $Y^j$, corresponding to the zero number $2j-1$, to parameters $\theta^j$, ${X^j}'$ and ${Y^j}'$, corresponding to the zero number $2j$. The mappings ${H^j}'$, $j=1,\ldots,k-1$ transfer triples $(\theta^j, {X^{j}}',{Y^j}')$ to $(t^{j+1},X^{j+1},Y^{j+1})$ which correspond to the zero number $2j+1$. The mapping $\widehat{I}_1$ transfers the triple $(\theta^k, {X^{k}}',{Y^k}')$ to the value of the corresponding solution for $t=T$.

If $\nu$ tends to infinity, $\delta$ tends to zero, the instants $t^j$ tend to the time instants of impacts of the solution of Eq. \eqref{eiwk2.1}, corresponding to the same initial conditions, uniformly with respect to $z_0$, $\widehat{I}_{0,1}\rightrightarrows I_{0,1}$ in the $C^1$ metrics. Due to the result of Lemma \ref{liwk10.2} and Eq. \eqref{eiwk10.6}, $H^j \rightrightarrows G^j$, ${H^j}'\rightrightarrows {\mathop{\rm id}\nolimits}$. It suffices to use Lemma \ref{liwk10.1} statement to finish the proof of Theorem \ref{thiwk2}. $\square$

\textbf{Remark.} Theorem \ref{thiwk2} applied to the results of Section 9 demonstrate that set of vibro-impact systems, satisfying Conditions 2 and 4---6, is non-empty and contains a subset open in ${\cal X}(J,n,T)$.

\bigskip

\noindent\textbf{11. Conclusion.}

\bigskip

We have considered vibro-impact systems in their general form. It was shown that the existence of an unstable periodic motion which passes near the delimiter without having an impact may imply chaos. The corresponding sufficient conditions can be written down explicitly. This shows that grazing is not a single bifurcation, but a combination of two bifurcations that can coexist. This phenomenon takes place both for impulse and "soft"\ models of impact dynamics.
\medskip

\noindent\textbf{Acknowledgements.} This work was supported by the UK Royal Society (joint research project of University of Aberdeen and Saint-Petersburg State University), by FCT research project PTDC/MAT/113470/2009, by Russian Foundation for Basic Researches, grant 12-01-00275-a and by the Chebyshev Laboratory (Department of Mathematics and Mechanics, Saint-Petersburg State University) under the grant 11.G34.31.0026 of the Government of the Russian Federation.

\end{document}